\documentclass[12pt,a4paper]{amsart}

\usepackage{t1enc}
\usepackage[latin1]{inputenc}
\usepackage[english]{babel}
\usepackage{hyperref}
\usepackage{graphicx}

\usepackage{latexsym}
\usepackage{amsmath,amsfonts,amssymb,amsthm}

\begin{document}

\title{On the Orbits of not Expansive Mappings in Metric Spaces}
\author{Sergio Venturini}
\address{
	S. Venturini:
	Dipartimento Di Matematica,
	Universit\`{a} di Bologna,
	\,\,Piazza di Porta S. Donato 5 ---I-40127 Bologna,
	Italy}
\email{sergio.venturini@unibo.it}


\def\R{{\rm I\kern-.185em R}}
\def\RR{\mathbb{R}}
\def\C{{\rm\kern.37em\vrule height1.4ex width.05em depth-.011em\kern-.37em C}}
\def\CC{\mathbb{C}}
\def\N{{\rm I\kern-.185em N}}
\def\NN{\mathbb{N}}
\def\Z{{\bf Z}}
\def\ZZ{\mathbb{Z}}
\def\Q{{\mathchoice{\hbox{\rm\kern.37em\vrule height1.4ex width.05em 
depth-.011em\kern-.37em Q}}{\hbox{\rm\kern.37em\vrule height1.4ex width.05em 
depth-.011em\kern-.37em Q}}{\hbox{\sevenrm\kern.37em\vrule height1.3ex 
width.05em depth-.02em\kern-.3em Q}}{\hbox{\sevenrm\kern.37em\vrule height1.3ex
width.05em depth-.02em\kern-.3em Q}}}}
\def\P{{\rm I\kern-.185em P}}
\def\H{{\rm I\kern-.185em H}}
\def\Aleph{\aleph_0}
\def\ALEPH#1{\aleph_{#1}}
\def\sset{\subset}\def\ssset{\sset\sset}
\def\bar#1{\overline{#1}}
\def\dim{\mathop{\rm dim}\nolimits}
\def\half{\textstyle{1\over2}}
\def\Half{\displaystyle{1\over2}}
\def\mlog{\mathop{\half\log}\nolimits}
\def\Mlog{\mathop{\Half\log}\nolimits}
\def\Det{\mathop{\rm Det}\nolimits}
\def\Hol{\mathop{\rm Hol}\nolimits}
\def\Aut{\mathop{\rm Aut}\nolimits}
\def\Re{\mathop{\rm Re}\nolimits}
\def\Im{\mathop{\rm Im}\nolimits}
\def\Ker{\mathop{\rm Ker}\nolimits}
\def\Fix{\mathop{\rm Fix}\nolimits}
\def\Exp{\mathop{\rm Exp}\nolimits}
\def\sp{\mathop{\rm sp}\nolimits}
\def\id{\mathop{\rm id}\nolimits}
\def\Rank{\mathop{\rm rk}\nolimits}
\def\Trace{\mathop{\rm Tr}\nolimits}
\def\Res{\mathop{\rm Res}\limits}
\def\cancel#1#2{\ooalign{$\hfil#1/\hfil$\crcr$#1#2$}}
\def\prevoid{\mathrel{\scriptstyle\bigcirc}}
\def\void{\mathord{\mathpalette\cancel{\mathrel{\scriptstyle\bigcirc}}}}
\def\n{{}|{}\!{}|{}\!{}|{}}
\def\abs#1{\left|#1\right|}
\def\norm#1{\left|\!\left|#1\right|\!\right|}
\def\nnorm#1{\left|\!\left|\!\left|#1\right|\!\right|\!\right|}
\def\upperint{\int^{{\displaystyle{}^*}}}
\def\lowerint{\int_{{\displaystyle{}_*}}}
\def\Upperint#1#2{\int_{#1}^{{\displaystyle{}^*}#2}}
\def\Lowerint#1#2{\int_{{\displaystyle{}_*}#1}^{#2}}
\def\rem #1::#2\par{\medbreak\noindent{\bf #1}\ #2\medbreak}
\def\proclaim #1::#2\par{\removelastskip\medskip\goodbreak{\bf#1:}
\ {\sl#2}\medskip\goodbreak}
\def\ass#1{{\rm(\rmnum#1)}}
\def\assertion #1:{\Acapo\llap{$(\rmnum#1)$}$\,$}
\def\Assertion #1:{\Acapo\llap{(#1)$\,$}}
\def\acapo{\hfill\break\noindent}
\def\Acapo{\hfill\break\indent}
\def\proof{\removelastskip\par\medskip\goodbreak\noindent{\it Proof.\/\ }}
\def\prova{\removelastskip\par\medskip\goodbreak
\noindent{\it Dimostrazione.\/\ }}
\def\risoluzione{\removelastskip\par\medskip\goodbreak
\noindent{\it Risoluzione.\/\ }}
\def\qed{{$\Box$}\par\smallskip}
\def\BeginItalic#1{\removelastskip\par\medskip\goodbreak
\noindent{\it #1.\/\ }}
\def\iff{if, and only if,\ }
\def\sse{se, e solo se,\ }
\def\rmnum#1{\romannumeral#1{}}
\def\Rmnum#1{\uppercase\expandafter{\romannumeral#1}{}}
\def\smallfrac#1/#2{\leavevmode\kern.1em
\raise.5ex\hbox{\the\scriptfont0 #1}\kern-.1em
/\kern-.15em\lower.25ex\hbox{\the\scriptfont0 #2}}
\def\Left#1{\left#1\left.}
\def\Right#1{\right.^{\llap{\sevenrm
\phantom{*}}}_{\llap{\sevenrm\phantom{*}}}\right#1}
\def\dimens{3em}
\def\symb[#1]{\noindent\rlap{[#1]}\hbox to \dimens{}\hangindent=\dimens}
\def\references{\bigskip\noindent{\bf References.}\bigskip}
\def\art #1 : #2 ; #3 ; #4 ; #5 ; #6. \par{#1, 
{\sl#2}, #3, {\bf#4}, (#5), #6.\par\smallskip}
\def\book #1 : #2 ; #3 ; #4. \par{#1, {\bf#2}, #3, #4.\par\smallskip}
\def\freeart #1 : #2 ; #3. \par{#1, {\sl#2}, #3.\par\smallskip}
\def\name{\hbox{Sergio Venturini}}
\def\snsaddress{\indent
\vbox{\bigskip\bigskip\bigskip
\name
\hbox{Scuola Normale Superiore}
\hbox{Piazza dei Cavalieri, 7}
\hbox{56126 Pisa (ITALY)}
\hbox{FAX 050/563513}}}
\def\cassinoaddress{\indent
\vbox{\bigskip\bigskip\bigskip
\name
\hbox{Universit\`a di Cassino}
\hbox{via Zamosch 43}
\hbox{03043 Cassino (FR)}
\hbox{ITALY}}}
\def\bolognaaddress{\indent
\vbox{\bigskip\bigskip\bigskip
\name
\hbox{Dipartimento di Matematica}
\hbox{Universit\`a di Bologna}
\hbox{Piazza di Porta S. Donato 5}
\hbox{40127 Bologna (BO)}
\hbox{ITALY}
\hbox{venturin@dm.unibo.it}
}}
\def\homeaddress{\indent
\vbox{\bigskip\bigskip\bigskip
\name
\hbox{via Garibaldi, 7}
\hbox{56124 Pisa (ITALY)}}}
\def\doubleaddress{
\vbox{
\hbox{\name}
\hbox{Universit\`a di Cassino}
\hbox{via Zamosch 43}
\hbox{03043 Cassino (FR)}
\hbox{ITALY}
\smallskip
\hbox{and}
\smallskip
\hbox{Scuola Normale Superiore}
\hbox{Piazza dei Cavalieri, 7}
\hbox{56126 Pisa (ITALY)}
\hbox{FAX 050/563513}}}
\def\sergio{{\rm\bigskip
\centerline{Sergio Venturini}
\leftline{\bolognaaddress}
\bigskip}}
\def\a{\alpha}
\def\bg{\beta}
\def\g{\gamma}
\def\G{\Gamma}
\def\dg{\delta}
\def\D{\Delta}
\def\e{\varepsilon}
\def\eps{\epsilon}
\def\z{\zeta}
\def\th{\theta}
\def\T{\Theta}
\def\k{\kappa}
\def\lg{\lambda}
\def\Lg{\Lambda}
\def\m{\mu}
\def\n{\nu}
\def\r{\rho}
\def\s{\sigma}
\def\Sg{\Sigma}
\def\ph{\varphi}
\def\Ph{\Phi}
\def\x{\xi}
\def\om{\omega}
\def\Om{\Omega}

\newtheorem{theorem}{Theorem}[section]
\newtheorem{proposition}{Proposition}[section]
\newtheorem{lemma}{Lemma}[section]
\newtheorem{corollary}{Corollary}[section]
\newtheorem{remark}{Remark}[section]
\newtheorem{definition}{Definition}[section]

\newtheorem{teorema}{Teorema}[section]
\newtheorem{proposizione}{Proposizione}[section]
\newtheorem{corollario}{Corollario}[section]
\newtheorem{osservazione}{Osservazione}[section]
\newtheorem{definizione}{Definizione}[section]
\newtheorem{esempio}{Esempio}[section]
\newtheorem{esercizio}{Esercizio}[section]
\newtheorem{congettura}{Congettura}[section]

\def\eeq{=}
\def\MSpace{X}
\def\MSpaceB{Y}
\def\MSubSpace{S}
\def\CSpace{X}
\def\CSpaceB{Y}
\def\LocSubSpace{W}
\def\KSpace{K}
\def\NSpace{Y}
\def\ngh{U}
\def\nghB{V}
\def\nghC{W}
\def\CompDivSet{E_{\pMap}}
\def\RecurrSet{S_{\pMap}}
\def\LimMapSet{\mathcal{G}_{\pMap}}
\def\LimMapSetEx{{\hat{\mathcal{G}}}_{\pMap}}
\def\distfunc{\delta}
\def\dist{\delta}
\def\radius{r}
\def\rad{\rho}
\def\ball#1#2{B_\dist\left({#1},{#2}\right)}
\def\uball#1#2{E_\dist\left({#1},{#2}\right)}
\def\Fam{\mathcal{F}}
\def\Ball{B}
\def\pMap{f}
\def\hMap{f}
\def\LMapA{u}
\def\LMapB{v}
\def\LMapC{w}
\def\pMapCA{h}
\def\ppa{x}
\def\ppb{y}
\def\ppc{z}
\def\ppd{w}
\def\ll{l}
\def\kk{k}
\def\nni{{n}}
\def\nnj{{m}}
\def\itn{\nu}
\def\itm{\mu}
\def\itx#1{{\kk_{#1}}}
\def\ity#1{{\ll_{#1}}}
\def\orbit{L}
\def\orbitb{K}
\def\limtoinfty#1{{\lim_{{#1}\to\infty}}}

\begin{abstract}
Let $\MSpace$ be a locally compact metric space and
let $\pMap:\MSpace\to\MSpace$ be a not expansive map.
We prove that for each $\ppa_0\in\MSpace$ the sequence
$\ppa_0,\pMap(\ppa_0),\pMap^2(\ppa_0),\ldots$
is either relatively compact in $\MSpace$ or compactly divergent in $\MSpace$.
As applications we study the structure of the functions which are limits
of the iterates of the map $\pMap$ and we prove the analyticity of the
set of $\pMap$-recurrent points when $\pMap:\MSpace\to\MSpace$ is a holomorphic
and $\MSpace$ is a complex hyperbolic spaces in the sense of Kobayashi.
\end{abstract}

\maketitle

\section{\label{section::Introduction}Introduction}
Let $\MSpace$ be a metric space with distance function $\distfunc_\MSpace$
and
let $\pMap:\MSpace\to\MSpace$ be a %continuous self-map.
\emph{not expansive} map,
that is a (necessarily continuous) map which satisfies
\begin{equation*}
\distfunc_\MSpace\bigl(\pMap(\ppa),\pMap(\ppb)\bigl)\leq\distfunc_\MSpace(\ppa,\ppb)
\end{equation*}
for each pairs of points $\ppa,\ppb\in\MSpace$.

The iterates of the map $\pMap$ are $\pMap^2=\pMap\circ\pMap$,
$\pMap^3=\pMap\circ\pMap\circ\pMap$
and so on.

For each $\ppa_0\in\MSpace$ the $\pMap$-\emph{orbit} of $\ppa_0$ is the sequence
\begin{equation*}
\ppa_0,\pMap(\ppa_0),\pMap^2(\ppa_0),\ldots
\end{equation*}

A sequence of points $\ppa_1, \ppa_2,\ldots$ in $\MSpace$ is said
\emph{compactly divergent} in $\MSpace$ if each compact subset $\KSpace\sset\MSpace$
the relation $\ppa_j\in\KSpace$ holds for a finite number of indexes $j$.

The main result of this paper is the following:

\begin{theorem}\label{thm::OutOut}
Let $\MSpace$ be a locally compact metric space %with countable basiss
and let $\pMap:\MSpace\to\MSpace$ be a not expansive map.

Then the $\pMap$-orbit of each point of $\MSpace$
is either relatively compact or compactly divergent in $\MSpace$.
\end{theorem}

Observe that we make no assumption on the completeness of $\MSpace$.

As an immediate consequence we obtain:

\begin{theorem}\label{thm::HyperbolicSpaces}
Let $\MSpace$ be a Kobayashi hyperbolic complex space
and let $\pMap:\MSpace\to\MSpace$ be a holomorphic map.

Then the $\pMap$-orbit of each point of $\MSpace$
is either relatively compact or compactly divergent in $\MSpace$.
\end{theorem}

For the definition of hyperbolicity in the sense of Kobayashi for a complex space
see, e.g.,  \cite{book:Kobayashi}, \cite{book:ComplexHyperbolicSpaces}
or section \ref{section::ComplexSpaces} below.

The paper is organized as follows.

Sections \ref{section::Preliminaries} and \ref{section::CalkaLemma} contains some
easy generalization of results already present in the literature that we
need for the proof of our Theorem \ref{thm::OutOut}, which is given in section
\ref{section::MainProof}.

In section \ref{section::LimitSet} we apply our main theorem 
to obtain a complete description of the structure of
the set of all the functions which are limit of iterates of a not expansive self-map
$\pMap:\MSpace\to\MSpace$, where $\MSpace$ is an
arbitrary locally compact metric space with countable basis.

Some further application to 
Kobayashi hyperbolic complex spaces are given 
in section \ref{section::ComplexSpaces}.

\section{\label{section::Preliminaries}Preliminaries}
In this paper we denote by
$\MSpace$ a metric space with distance function $\distfunc_\MSpace$
and $\pMap:\MSpace\to\MSpace$ will be
a not expansive map of $\MSpace$ in itself.

For each $\ppa\in\MSpace$ and each $\radius>0$ we denote by $\Ball_\MSpace(\ppa,\radius)$
the open ball in $\MSpace$ of center $\ppa$ and radius $\radius$
and
for each subset $\KSpace\subset\MSpace$ % and each $\radius>0$
we set
\begin{equation*}
\KSpace_\radius=\bigcup_{\ppc\in\KSpace}\Ball_\MSpace(\ppc,\radius)
\end{equation*}
that is $\ppa\in\KSpace_\radius$ if, and only if,
$\distfunc_\MSpace(\ppa,\ppc)<\radius$ for some $\ppc\in\KSpace$.

Let us begin with the following simple observation.
\begin{proposition}\label{thm::ConvLemma}
Let $\MSpace$ and $\NSpace$ be two metric space.
Let $\pMap_\nni:\MSpace\to\NSpace$ be a sequence of not expansive mappings
and let $\ppa\in\MSpace$.

If for some sequence $\ppa^0_\nni$ converging to $\ppa$ we have
\begin{equation*}
\limtoinfty\nni\pMap_\nni(\ppa^0_\nni)=\ppc
\end{equation*}
with $\ppc\in\NSpace$ then for each sequence $\ppa_\nni$ converging to $\ppa$ we also have
\begin{equation*}
\limtoinfty\nni\pMap_\nni(\ppa_\nni)=\ppc.
\end{equation*}
In particular we have
\begin{equation*}
\limtoinfty\nni\pMap_\nni(\ppa)=\ppc.
\end{equation*}
\end{proposition}

\proof
Let denote by $\distfunc_\MSpace$ and $\distfunc_\NSpace$
the distance functions respectively on $\MSpace$ and $\NSpace$.
Then
\begin{eqnarray*}
\distfunc_\NSpace\bigl(\ppc,\pMap_\nni(\ppa_\nni)\bigr)&\leq&
\distfunc_\NSpace\bigl(\ppc,\pMap_\nni(\ppa^0_\nni)\bigr)
	+\distfunc_\NSpace\bigl(\pMap_\nni(\ppa^0_\nni),\pMap_\nni(\ppa_\nni)\bigr)
\\
&\leq&
\distfunc_\NSpace\bigl(\ppc,\pMap_\nni(\ppa^0_\nni)\bigr)
	+\distfunc_\MSpace(\ppa^0_\nni,\ppa_\nni)
\\
&\leq&
\distfunc_\NSpace\bigl(\ppc,\pMap_\nni(\ppa^0_\nni)\bigr)
	+\distfunc_\MSpace(\ppa^0_\nni,\ppa)
	+\distfunc_\MSpace(\ppa,\ppa_\nni).
\end{eqnarray*}
Taking the limit as $n\to\infty$, observing that
\begin{equation*}
\limtoinfty\nni\distfunc_\NSpace\bigl(\ppc,\pMap_\nni(\ppa^0_\nni)\bigr)
	=\limtoinfty\nni\distfunc_\MSpace(\ppa^0_\nni,\ppa)
	=\limtoinfty\nni\distfunc_\MSpace(\ppa,\ppa_\nni)=0,
\end{equation*}
we obtain
\begin{equation*}
\limtoinfty\nni\distfunc_\NSpace\bigl(\ppc,\pMap_\nni(\ppa_\nni)\bigr)=0,
\end{equation*}
as desired.

\qed

In particular we obtain
\begin{lemma}\label{thm::CompToRecurrent}
Let $\MSpace$ be a metric space and let
$\pMap:\MSpace\to\MSpace$ be a not expansive
self-map.
Let $\ppa,\ppb\in\MSpace$ and let $\itx\itn$ be a sequence of positive integers
such that $\itx\itn$ and $\itx{\itn+1}-\itx\itn$ are both increasing sequences.
Then
\begin{equation*}
\limtoinfty\itn\pMap^{\itx\itn}(\ppa)=\ppb
\Longrightarrow
\limtoinfty\nu\pMap^{\itx{\itn+1}-\itx\itn}(\ppb)=\ppb.
\end{equation*}
In particular it follows that $\ppb\in\MSpace$ is $\pMap$-recurrent.
\end{lemma}

\proof
Indeed, by hypotheses,
\begin{eqnarray*}
&&\limtoinfty\nu\pMap^{\itx\itn}(\ppa)=\ppb,\\
&&\limtoinfty\nu\pMap^{\itx{\itn+1}-\itx\itn}\bigl(\pMap^{\itx\itn}(\ppa)\bigr)=
\limtoinfty\nu\pMap^{\itx{\itn+1}}(\ppa)=\ppb
\end{eqnarray*}
and hence, by the previous proposition,
\begin{equation*}
\limtoinfty\nu\pMap^{\itx{\itn+1}-\itx\itn}(\ppb)=\ppb,
\end{equation*}
as desired.

\qed

We shall need of a topological version of the Ascoli-Arzela theorem.
Let $\MSpace$ and $\MSpaceB$ two topological spaces.
Let us recall that 
a family $\mathcal{F}\subset C(\MSpace,\MSpaceB)$
is \emph{evenly continuous} if
for every $\ppa\in\MSpace$, $\ppb\in\MSpaceB$,
and every neighbourhood $\nghB$ of $\ppb$ in $\MSpaceB$
there are a neighbourhood $\ngh$ of $\ppa$ in $\MSpace$
and a neighbourhood $\nghC$ of $\ppb$ in $\MSpaceB$
such that for every $\pMap\in\mathcal{F}$
\begin{equation*}
\pMap(\ppa)\in\nghC\Longrightarrow\pMap(\ngh)\subset\nghB.
\end{equation*}

Then the topological Ascoli-Arzela theorem
given in \cite[7.21]{book:KelleyGenTopol} is
\begin{theorem}\label{thm:AACompactness}
Let $\MSpace$ be a regular locally compact topological space
and $\MSpaceB$ a regular Hausdorff topological space.
Then a family $\mathcal{F}\subset C(\MSpace,\MSpaceB)$ is
relatively compact in $C(\MSpace,\MSpaceB)$ if, and only if,
it is evenly continuous and 
$$\bigl\{\pMap(\ppa)\mid\pMap\in\mathcal{F}\bigr\}$$
is relatively compact in $\MSpaceB$ for all $\ppa\in\MSpace$.

In particular if $\MSpaceB$ is compact
then $\mathcal{F}\subset C(\MSpace,\MSpaceB)$
relatively compact in $C(\MSpace,\MSpaceB)$ if, and only if,
it is evenly continuous.
\end{theorem}

The following results are straightforward generalization of some
results due to Loeb and Vigu\'e (\cite{article:LoebVig}).

\begin{proposition}
Let $\MSpace$ and $\NSpace$ be two metric space with distance function
respectively $\distfunc_\MSpace$ and $\distfunc_\NSpace$.

Let $\pMap:\MSpace\to\NSpace$ be a continuos map.
Assume that the image $\pMap(\MSpace)$ is dense in $\NSpace$ and
that
\begin{equation*}
\distfunc_\MSpace(\ppa,\ppb)\leq\distfunc_\NSpace\bigl(\pMap(\ppa),\pMap(\ppb)\bigr)
\end{equation*}
for each pair of point $\ppa,\ppb\in\MSpace$.

Let $\ppa_0\in\MSpace$  and let $\radius>0$ be given.
If $\bar{\Ball_\MSpace(\ppa_0,\radius)}$ is complete (as metric space) then
\begin{equation*}
\Ball_\NSpace\bigl(\pMap(\ppa_0),\radius\bigr)\sset\pMap\bigl(\Ball_\MSpace(\ppa_0,\radius)\bigr).
\end{equation*}

\end{proposition}

\proof
Let $\ppb\in\Ball_\NSpace\bigl(\pMap(\ppa_0),r\bigr)$.
We need to prove that there exits $\ppa\in\Ball_\MSpace(\ppa_0,\radius)$
such that $\pMap(x)=y$.

Since $\pMap(\MSpace)$ is dense in $\NSpace$ there exist a sequence $\ppa_\nni\in\MSpace$
such that $\pMap(\ppa_\nni)\to\ppb$. It is not restrictive to assume that
$\pMap(\ppa_\nni)\in\Ball_\NSpace\bigl(\pMap(\ppa_0),\radius\bigr)$.
We then have
\begin{equation*}
\distfunc_\MSpace(\ppa_\nni,\ppa_0)\leq
\distfunc_\NSpace\bigl(\pMap(\ppa_\nni),\pMap(\ppa_0)\bigr)<\radius,
\end{equation*}
that is $\ppa_\nni\in\Ball_\MSpace(\ppa_0,\radius)$.
We also have
\begin{equation*}
\distfunc_\MSpace(\ppa_\nni,\ppa_\nnj)\leq
\distfunc_\NSpace\bigl(\pMap(\ppa_\nni),\pMap(\ppa_\nnj)\bigr)
\leq
\distfunc_\NSpace\bigl(\pMap(\ppa_\nni),\ppb\bigr)+
\distfunc_\NSpace\bigl(\ppb,\pMap(\ppa_\nnj)\bigr)
\end{equation*}
and therefore the sequence $\ppa_\nni$ is a Cauchy sequence in $\Ball_\MSpace(\ppa_0,\radius)$.

Since $\bar{\Ball_\MSpace(\ppa_0,\radius)}$ is complete there exist
$\ppa\in\bar{\Ball_\MSpace(\ppa_0,\radius)}$ such that the sequence $\ppa_n$ converges to $\ppa$.
But
\begin{equation*}
\distfunc_\MSpace(\ppa,\ppa_0)=\limtoinfty\nni\distfunc_\MSpace(\ppa_\nni,\ppa_0)
\end{equation*}
and
\begin{equation*}
\limtoinfty\nni\distfunc_\MSpace(\ppa_\nni,\ppa_0)
\leq
\limtoinfty\nni\distfunc_\NSpace\bigl(\pMap(\ppa_\nni),\pMap(\ppa_0)\bigr)=
\distfunc_\NSpace\bigl(\ppb,\pMap(\ppa_0)\bigr)<\radius,
\end{equation*}

that is $\ppa\in\Ball_\MSpace(\ppa_0,\radius)$ and since $\pMap$ is continuous
\begin{equation*}
\pMap(\ppa)=\pMap\left(\limtoinfty\nni\ppa_\nni\right)=
\limtoinfty\nni\pMap(\ppa_\nni)=\ppb,
\end{equation*}
as desired.

\qed

We say that a metric space $\MSpace$ with distance function $\distfunc_\MSpace$
is \emph{locally complete} if for each $\ppa\in\MSpace$ there exists $\radius>0$ such that
$\bar{\Ball_\MSpace(\ppa,\radius)}$ is complete (as metric space).

Of course each locally compact metric space is locally complete.

\begin{proposition}\label{thm::MonoToIso}
Let $\MSpace$ be a locally complete metric space
and let $\pMap:\MSpace\to\MSpace$ be a not expansive map.
Assume that for an increasing sequence of positive integer $\itx\itn$
the sequence $\pMap^{\itx\itn}$ converges pointwise to the identity map of $\MSpace$.
Then $\pMap$ is a surjective isometry.
\end{proposition}

\proof
Let $\itx\itn$ be an increasing sequence of positive integers such that
the sequence $\pMap^{\itx\itn}$ converges to the identity map of $\MSpace$.

Let $\ppa,\ppb\in\MSpace$.
Then
\begin{eqnarray*}
\distfunc_\MSpace(\ppa,\ppb)&=&
\limtoinfty\itn\distfunc_\MSpace\bigl(\pMap^{\itx\itn}(\ppa),\pMap^{\itx\itn}(\ppb)\bigr)
\\
&=&\limtoinfty\itn
\distfunc_\MSpace\Bigl(\pMap^{\itx\itn-1}\bigl(\pMap(\ppa)\bigr),\pMap^{\itx\itn-1}\bigl(\pMap(\ppb\bigr)\Bigr)
\\
&\leq&
\distfunc_\MSpace\bigl(\pMap(\ppa),\pMap(\ppb)\bigr).
\end{eqnarray*}
Since we also have $\distfunc_\MSpace(\ppa,\ppb)\geq\distfunc_\MSpace\bigl(\pMap(\ppa),\pMap(\ppb)\bigr)$
then it follows that
\begin{equation*}
\distfunc_\MSpace(\ppa,\ppb)=\distfunc_\MSpace\bigl(\pMap(\ppa),\pMap(\ppb)\bigr),
\end{equation*}
that is the map $\pMap$ is an isometry.
By induction on $\kk$ it follows that $\pMap^\kk$ is an isometry too.

It remains to show that the map $\pMap$ is surjective.

For each $\ppa\in\MSpace$ and each $\itn>0$ we have
\begin{equation*}
0=\limtoinfty\itm\distfunc_\MSpace\bigl(\ppa,\pMap^{\itx\itm}(\ppa)\bigr)=
\limtoinfty\itm\distfunc_\MSpace\Bigl(\ppa,\pMap^{\itx\itn}\bigl(\pMap^{\itx\itm-\itx\itn}(\ppa)\bigr)\Bigr).
\end{equation*}
It follows that for each $\itn$ the image $\pMap^{\itx\itn}(\MSpace)$ is dense in $\MSpace$.

Let now $\ppa\in\MSpace$ be arbitrary. Choose $\radius>0$ such that $\bar{\Ball_\MSpace(\ppa,\radius)}$
is complete. Then, by the previous proposition for each $\itn>0$ 
\begin{equation*}
\Ball_\MSpace\bigl(\pMap^{\itx\itn}(\ppa),\radius\bigr)\sset\pMap^{\itx\itn}\bigl(\Ball_\MSpace(\ppa,\radius)\bigr)
\sset\pMap(\MSpace).
\end{equation*}
But for $\itn>0$ large enough we have $\ppa\in\Ball_\MSpace\bigl(\pMap^{\itx\itn}(\ppa),\radius\bigr)$
and hence $\ppa\in\pMap(\MSpace)$.
Since $\ppa\in\MSpace$ is arbitrary it follows that the map $\pMap$ is surjective.

\qed

\section{\label{section::CalkaLemma}A lemma of Ca\l ka}
{ %begin unit Calka
The main result in this section (Theorem \ref{thm::CalkaLemma})
is a reformulation a results given by Ca\l ka
in \cite[Lemma 3.1 pag. 222]{article:CalkaBoundedOrbits}.

Let $\dist:\NN\times\NN\to[0,+\infty[$ be a distance function on $\NN$.
For each $n\in\NN$ and each $\rad>0$ we set
\begin{eqnarray*}
\ball n\rad&=&\bigl\{k\in\NN\mid \dist(k,n)<\rad\bigr\},\\
\uball n\rad&=&\bigcup_{k=0}^n\ball n\rad.
\end{eqnarray*}

Of course $\ball n\rad\subset\uball n\rad$
and
\begin{equation*}
\uball n\rad\subset\uball {n+1}\rad\subset\uball {n+2}\rad\subset\cdots.
\end{equation*}

\begin{theorem}[Ca\l ka lemma]\label{thm::CalkaLemma}
Let $\dist:\NN\times\NN\to[0,+\infty[$ be a distance function on $\NN$ such that
\begin{equation*} %\label{eq:distanceProp}
\dist(n+1,m+1)\geq\dist(n,m)
\end{equation*}
for each $n,m\in\NN$.

Assume that for some $N\in\NN$ and $\rad>0$ the ball  $\ball 0\rad$ is infinite
and
\begin{equation*}
\ball 0\rad\subset\uball N{\rad/2}.
\end{equation*}
Then 
\begin{equation*}
\NN=\uball M\rad
\end{equation*}
for some $M\in\NN$.
\end{theorem}

For the proof we need of the following lemma.

\begin{lemma}\label{lemma:subCalkaLemma}
Let $\dist:\NN\times\NN\to[0,+\infty[$ be a distance function on $\NN$
as in Theorem \ref{thm::CalkaLemma}.
Let $n,\nu,m\in\NN$ and $\rad>0$ satisfying
\begin{equation*}
n<\nu<m,\ \nu\not\in\uball n\rad,\ m\in\ball n\rad.
\end{equation*}
Then 
\begin{eqnarray*}
&&\nu<m-n,\\
&&\dist(m-n,0)\leq\dist(m,n)<\rad.
\end{eqnarray*}
\end{lemma}

\proof
Since $\dist(n+1,m+1)\geq\dist(n,m)$ for each $n,m\in\NN$ it follows that
the sequence
\begin{equation*}
	j\mapsto\dist(m-n+j,j)
\end{equation*}
is not decreasing and hence
\begin{equation*}
	\dist(m-n,0)\leq\dist(m,n)<\rad.
\end{equation*}

we also have
\begin{equation*}
0\leq j\leq n\ \Longrightarrow\ \dist(m-n+j,j)\leq\dist(m,n)<\rad,
\end{equation*}
and hence, when $0\leq j\leq n$,
\begin{equation*}
m-n+j\in\ball j\rad\subset\uball n\rad,
\end{equation*}
that is
\begin{equation*}
m-n\leq k\leq m\ \Longrightarrow\ k\in\uball n\rad.
\end{equation*}
Being $\nu<m$ and also $\nu\not\in\uball M\rad$
necessarily $\nu<m-n$, as required.

\qed

Let now $N\in\NN$ and $\rad>0$ such that
the ball $\ball 0\rad$ is infinite
and
$ %\begin{equation*}
\ball 0\rad\subset\uball N{\rad/2},
$ %\end{equation*}
that is
\begin{equation*}
\ball0\rad\subset\bigcup_{k=0}^n\ball n{/2}.
\end{equation*}
As $\ball0\rad$ contains infinite positive integers
it follows that $\ball{n_0}{\rad/2}$ also contains infinite positive integers
for some $n_0\leq N$.

Observe that if  $k\in\ball{n_0}{\rad/2}$ and $k\geq n_0$
then the sequence $j\mapsto\dist(k-n_0+j,j)$ is not increasing and hence
\begin{equation*}
\dist(k-n_0,0)\geq\dist(k,n_0)<\rad/2,
\end{equation*}
that is the $\ball0{\rad/2}$ contains all the infinite positive integers
$k-n_0$ with $k\in\ball{n_0}{\rad/2}$ and $k\geq n_0$.

Let now $M\in\NN$ with $M>N$ and $\dist(0,M)<\rad/2$.
We end the proof of Theorem \ref{thm::CalkaLemma}
showing that $\NN=\uball M\rad$.

Assume by contradiction that exists $\nu\in\NN$ such that
$\nu\not\in\uball M\rad$.
Clearly $\nu>M$.

We have already observed that the ball $\ball0{\rad/2}$ is infinite,
so let $m_0$ be the first positive integer which satisfies
$m_0>\nu$ and $m_0\in\ball0{\rad/2}$.

Then $M<\nu<m_0$ and the triangle inequality implies
\begin{equation*}
\dist(m_0,M)\leq \dist(m_0,0)+\dist(0,M)<\rad/2+\rad/2=\rad.
\end{equation*}
Lemma \ref{lemma:subCalkaLemma} implies that $\nu<m_0-M$ and
\begin{equation*}
\dist(m_0-M,0)\leq\dist(m_0,M)<\rad,
\end{equation*}
that is $m_0-M\in\ball0\rad$.

By our hypotheses $\ball0\rad\subset\uball N{\rad/2}$ and hence 
$m_0-M\in\ball{n}{/2}$ for some $n\in\NN$ satisfying $n\leq N$.

Clearly $n\leq N<M<\nu<m_0-M$ and
\begin{equation*}
\dist(m_0-M,n)<\rad/2<\rad.
\end{equation*}
We apply lemma \ref{lemma:subCalkaLemma} again and obtain that
setting $m_1=m_0-M-n$ then $\nu<m_1$ and
\begin{equation*}
\dist(m_1,0)=\dist(m_0-M-n,0)=\dist(m_0-M,n)<\rad/2
\end{equation*}
and this contradict the choice of $m_0$ as the smallest positive integer
which satisfies $m>\nu$ and $\dist(m,0)<\rad/2$.

} %end unit Calka
\section{\label{section::MainProof}Proof of the main theorem}
{ %begin unit
\def\pBase{{\ppa_0}}
\def\PB{E}

Let us begin with the following particular case
of Theorem \ref{thm::OutOut}.

\begin{lemma}
Let $\MSpace$ be a locally compact metric space
and let $\pMap:\MSpace\to\MSpace$ be a surjective isometry.
Then the $\pMap$-orbit of each $\pMap$-recurrent point is
relatively compact in $\MSpace$.
\end{lemma}

{
\proof
Let $\pBase\in\MSpace$ be a $\pMap$-recurrent point of $\MSpace$
and let $\orbit$ be the $\pMap$-orbit of the point $\pBase$.

If the map $\NN\ni\nni\mapsto\pMap^\nni(\ppa)\in\MSpace$
is not injective $\orbit$ is finite and hence compact.

Assume hence that $\pMap^{\nni}(\pBase)\neq\pMap^{\nnj}(\pBase)$
when $\nni\neq\nnj$.
Consider the distance function on $\NN$ defined by the formula
\begin{equation*}
\dist(\nni,\nnj)=\distfunc_\MSpace\bigl(\pMap^\nni(\pBase),\pMap^\nnj(\pBase)\bigr)
\end{equation*}
and choose $\rad>0$ in such a way that the ball $\Ball_\MSpace	(\pBase,\rad)$
is relatively compact in $\MSpace$.

Since $\pMap$ is an isometry we have
\begin{equation*}
\dist(\nni+1,\nnj+1)=\dist(\nni,\nnj)
\end{equation*}
for each $\nni,\nnj\in\NN$.

Let define $\ball n\rad$ and $\uball n\rad$ as in the previous section.

Since the point $\pBase$ is $\pMap$-recurrent the ball $\ball 0\rad$
is infinite.

Let $$\PB=\bigl\{\pMap^\nni(\pBase)\mid\nni\in\ball 0\rad\bigr\}$$
and let $\bar{\PB}$ be the closure of $\PB$ in $\MSpace$.
Clearly we have
\begin{equation*}
\bar\PB\subset\bigcup_{\nni\in\ball 0\rad}\Ball_{\MSpace}(\pMap^\nni(\pBase),\rad/2)
\subset\bigcup_{\nni\in\NN}\Ball_{\MSpace}(\pMap^\nni(\pBase),\rad/2).
\end{equation*}
Since $\bar{\PB}\subset\Ball_\MSpace(\pBase,\rad)$ it follows that
$\bar{\PB}$ is compact and hence there exists $N\in\NN$ such that
\begin{equation*}
\PB\subset\bar\PB\subset\bigcup_{\nni=0}^N\Ball_{\MSpace}(\pMap^\nni(\pBase),\rad/2),
\end{equation*}
and hence
\begin{equation*}
\ball 0\rad\subset\uball N{\rad/2}.
\end{equation*}
Theorem \ref{thm::CalkaLemma} implies that for some $M>0$
\begin{equation*}
\NN\subset\uball M{\rad},
\end{equation*}
that is
\begin{equation*}
\orbit\subset\bigcup_{\nni\in\ball 0\rad}^M\Ball_{\MSpace}(\pMap^\nni(\pBase),\rad).
\end{equation*}
But the map $\pMap$ is a surjective isometry and hence for each $\nni\in\NN$
\begin{equation*}
\Ball_{\MSpace}\bigl(\pMap^\nni(\pBase),\rad\bigr)=\pMap^\nni\bigl(\Ball_{\MSpace}(\pBase,\rad)\bigr)
\end{equation*}
is relatively compact.
It follows that  the orbit $\orbit$ is contained in a finite union of relatively compact subset of
$\MSpace$ and hence is a relatively compact subset of $\MSpace$, as required.

}
\qed

\begin{lemma}
Let $\MSpace$ be a locally compact metric space
and let $\pMap:\MSpace\to\MSpace$ be a not expansive map.
Then the $\pMap$-orbit of each $\pMap$-recurrent point is
relatively compact in $\MSpace$.
\end{lemma}

{
\def\subspace{E}
\proof
Let $\pBase\in\MSpace$ be a $\pMap$-recurrent point
and let $\orbit\subset\MSpace$ be its $\pMap$-orbit.
By definition of $\pMap$-recurrent point there exist
an increasing sequence of positive integers $\itx\itn$ such that
\begin{equation*}
\limtoinfty\itn\pMap^{\itx\itn}(\pBase)=\pBase.
\end{equation*}
Define
\begin{equation*}
\subspace=\bigl\{\ppa\in\MSpace\mid\limtoinfty\itn\pMap^{\itx\itn}(\ppa)=\ppa\bigr\}.
\end{equation*}
Then $\subspace\neq\void$ because $\pBase\in\subspace$.

If $\ppa\in\subspace$ then
\begin{equation*}
\limtoinfty\itn\pMap^{\itx\itn}\bigl(\pMap(\ppa)\bigr)=
\limtoinfty\itn\pMap\bigl(\pMap^{\itx\itn}(\ppa)\bigr)=
\pMap\bigl(\limtoinfty\itn\pMap^{\itx\itn}(\ppa)\bigr)=
\pMap(\ppa),
\end{equation*}
and hence $\pMap(\ppa)\in\subspace$, that is $\pMap(\subspace)\subset\subspace$
and $\orbit\subset\subspace$.

We claim that $\subspace$ is closed in $\MSpace$.
Indeed let $\ppa\in\bar{\subspace}$ and let $\e>0$.
Choose $\ppb\in\subspace$ which satisfies $\distfunc_\MSpace(\ppa,\ppb)<\e$.
If $\itn\in\NN$ is large enougth we have
$\distfunc_\MSpace\bigl(\pMap^{\itx\itn}(\ppb),\ppb\bigr)<\e$
and
\begin{eqnarray*}
	\distfunc_\MSpace\bigl(\pMap^{\itx\itn}(\ppa),\ppa\bigr)&\leq&
		\distfunc_\MSpace\bigl(\pMap^{\itx\itn}(\ppa),\pMap^{\itx\itn}(\ppb)\bigr)+
		\distfunc_\MSpace\bigl(\pMap^{\itx\itn}(\ppb),\ppa\bigr)\\
	&\leq&
		\distfunc_\MSpace(\ppa,\ppb)+
		\distfunc_\MSpace\bigl(\ppa,\pMap^{\itx\itn}(\ppb)\bigr)\\
	&\leq&
		\distfunc_\MSpace(\ppa,\ppb)+
		\distfunc_\MSpace(\ppa,\ppb)+
		\distfunc_\MSpace\bigl(\ppb,\pMap^{\itx\itn}(\ppb)\bigr)\\
	&\leq&3\,\e.
\end{eqnarray*}
Since $\e>0$ can be chosen arbitrarily small it follows that
\begin{equation*}
\limtoinfty\itn\pMap^{\itx\itn}(\ppa)=\ppa,
\end{equation*}
that is $\ppa\in\subspace$.

Thus $\subspace$ is a locally compact because is a closed subset of
the locally compact space $\MSpace$.

We have $\pMap(\subspace)\subset\subspace$ and the sequence $\pMap^{\itx\itn}$
converges pointwise to the identity map of $\subspace$. 
Proposition \ref{thm::MonoToIso} implies that $\pMap$ is a surjective isometry
and the previous lemma yields that the orbit $\orbit$ of the $\pMap$-recurrent point
$\pBase$ is relatively compact in $\subspace$ and hence also in $\MSpace$.

\qed }

We are now able to prove Theorem \ref{thm::OutOut}.
Let $\ppa$ be a point of $\MSpace$.
Assume that the $\pMap$-orbit of $\ppa$ is not compactly divergent.
Then there exist $\ppb\in\MSpace$ and an increasing sequence
of positive integers $\itx\itn$ such that 
\begin{equation*}
\limtoinfty\itn\pMap^{\itx\itn}(\ppa)=\ppb.
\end{equation*}
It is not restrictive to assume that also the sequence
$\itx{\itn+1}-\itx\itn$ is increasing.
Then Lemma \ref{thm::CompToRecurrent} implies that the point $\ppb$ is $\pMap$-recurrent.
Let $\orbitb$ be the $\pMap$-orbit of $\ppb$. By the previous lemma $\orbitb$ is
relatively compact in $\MSpace$.

Let $\e>0$ be small enough in such a way that $\orbitb_\e$ %the set
is relatively compact in $\MSpace$.

Choose $\itn_0\in\NN$ which satisfies
$\distfunc_\MSpace(\pMap^{\itx{\itn_0}}(\ppa),\ppb)<\e$.
Then for each $\itn>\itn_0$ we have
$\pMap^{\itx\itn-\itx{\itn_0}}(\ppb)\in\orbitb$ and
\begin{eqnarray*}
\distfunc_\MSpace\bigl(\pMap^{\itx\itn}(\ppa),\pMap^{\itx\itn-\itx{\itn_0}}(\ppb)\bigr)
&\leq&
	\distfunc_\MSpace\bigl(\pMap^{\itx\itn-\itx{\itn_0}}\bigl(\pMap^{\itx{\itn_0}}(\ppa)\bigr),
		\pMap^{\itx\itn-\itx{\itn_0}}(\ppb)\bigr)\\
&\leq&
	\distfunc_\MSpace\bigl(\pMap^{\itx{\itn_0}}(\ppa),(\ppb)\bigr)<\e,
\end{eqnarray*}
that is
$\pMap^{\itx\itn}(\ppa)\in\orbitb_\e$.
It follows that the orbit $\orbit$ is relatively compact in $\MSpace$ being 
contained in 
\begin{equation*}
\orbitb_\e\cup\{\ppa\}\cup\{\pMap(\ppa)\}\cup\cdots\cup\{\pMap^{\itn_0}(\ppa)\},
\end{equation*}
which is clearly a relatively compact subset of $\MSpace$.

} %end unit

\section{\label{section::LimitSet}Limits of iterates}
In this section we assume that $\MSpace$ is
a locally compact metric space with countable basis.

The main result of this section, Theorem \ref{thm::GroupStructBase},
is a complete description of the structure of the maps
which are limit of sequences of iterates of a not expansive map
$\pMap:\MSpace\to\MSpace$ and 
is inspired to the results of Abate on the study of the limit points
of the iterates of an holomorphic map on taut complex manifolds:
see \cite[Theorem 2.1.29 pag. 143]{book:AbateBook}.

We denote by $\hat{\MSpace}=\MSpace\cup\{\infty\}$
the Alexandroff compactification of
the locally compact but not compact space $\MSpace$
and we set $\hat{\MSpace}=\MSpace$ if $\MSpace$ is compact.

If $\NSpace$ is an other metric space we denote by $C(\MSpace,\NSpace)$ the
set of all the continuous maps from $\MSpace$ to $\NSpace$ endowed with the
compact open topology.
Then it is straightforward to prove that the composition map
\begin{equation*}
C(\MSpace,\NSpace)\times C(\MSpace,\NSpace)
\ni(\LMapA,\LMapB)\mapsto\LMapA\circ\LMapB\in
C(\MSpace,\NSpace)
\end{equation*}
is continuous.

The proof of the theorem above follows the same lines of
\cite[Lemma 1.2 pag. 791]{article:AbateCharHypMan}.

\begin{proposition}\label{thm::AbateComp}
Let $\MSpace$ and $\NSpace$ be two locally compact metric space with countable base.
Then the family $\mathcal{F}\subset C(\MSpace,\NSpace)$ 
of all not expansive maps from $\MSpace$ to $\NSpace$
is relatively compact in $C(\MSpace,\hat\NSpace)$.
\end{proposition}

\proof
Since $\hat\NSpace$ is compact,
by Theorem \ref{thm:AACompactness}
it suffices to prove that the family $\mathcal{F}$ 
is evenly continuous in $C(\MSpace,\hat\NSpace)$.

Let denote by $\distfunc_\MSpace$ and $\distfunc_\MSpaceB$
the distance functions respectively on $\MSpace$ and $\MSpaceB$.

Let $\ppa\in\MSpace$, $\ppb\in\hat\MSpaceB$ and
$\nghB$ a neighbourhood of $\ppb$ in $\MSpaceB$.

Suppose first that $\ppa\neq\infty$, that is $\ppa\in\MSpaceB$.
Then choose $\rad>0$ small enough satisfying
$\Ball_{\MSpaceB}(\ppb,2\rad)\subset\nghB$
and set $\ngh=\Ball_{\MSpace}(\ppa,\rad)$ and
$\nghC=\Ball_{\MSpaceB}(\ppb,\rad)$.

Let $\pMap\in\Fam$ and suppose that $\pMap(\ppa)\in\nghC$,
that is $\distfunc_\MSpaceB\bigl(\pMap(\ppa),\ppb\bigr)<\rad$.
If $\ppc\in\ngh$ then $\distfunc_\MSpaceB(\ppc,\ppa)<\rad$ and
\begin{equation*}
\distfunc_\MSpaceB\bigl(\pMap(\ppc),\ppb\bigr)\leq
\distfunc_\MSpaceB\bigl(\pMap(\ppc),\pMap(\ppa)\bigr)+
\distfunc_\MSpaceB\bigl(\pMap(\ppa),\ppb\bigr)
<\distfunc_\MSpaceB(\ppc,\ppa)+
\rad<2\,\rad,
\end{equation*}
that is $\pMap(\ppc)\in\nghB$.
Since $\ppc\in\ngh$ is arbitrary we have $\pMap(\ngh)\subset\nghB$.

If $\MSpaceB$ is compact we are done.
Assume hence that $\MSpaceB$ is not compact and  $\ppb=\infty$.
Let $\KSpace\subset\MSpaceB$ a compact set such that
$\MSpaceB\setminus\KSpace\subset\nghB$.

Choose $\rad>0$ in such a way that $\KSpace_{\rad}$
is relatively compact in $\MSpaceB$ and set
$\ngh=\Ball_{\MSpace}(\ppa,\rad)$ and
$\nghC=\MSpaceB\setminus\bar{\KSpace_{\rad}}\cup\{\infty\}$.

Let $\pMap\in\Fam$ and suppose that $\pMap(\ppa)\in\nghC$.
If $\ppb\in\KSpace$ and $\ppc\in\ngh$ then
$\distfunc_\MSpaceB\bigl(\pMap(\ppa),\ppb\bigr)\geq\rad$,
$\distfunc_\MSpaceB(\ppa,\ppc)<\rad$
and
\begin{equation*}
\distfunc_\MSpaceB\bigl(\pMap(\ppa),\ppb\bigr)\leq
\distfunc_\MSpaceB\bigl(\pMap(\ppa),\pMap(\ppc)\bigr)+
\distfunc_\MSpaceB\bigl(\pMap(\ppc),\ppb\bigr)\leq
\distfunc_\MSpaceB(\ppa,\ppc)+
\distfunc_\MSpaceB\bigl(\pMap(\ppc),\ppb\bigr)
\end{equation*}
and hence
\begin{equation*}
\distfunc_\MSpaceB\bigl(\pMap(\ppc),\ppb\bigr)\geq
\distfunc_\MSpaceB\bigl(\pMap(\ppa),\ppb\bigr)
-\distfunc_\MSpaceB(\ppa,\ppc)>\rad-\rad=0,
\end{equation*}
that is $\pMap(\ppc)\neq\ppb$.

Since $\ppb\in\KSpace$ and $\ppc\in\ngh$
are arbitrary then
$\pMap(\ngh)\cap\KSpace=\void$ and hence
$\pMap(\ngh)\subset\MSpaceB\setminus\KSpace\subset\nghB$.

\qed

It is an immediate consequence of the theorem above
that the topology of the pointwise convergence and
the compact open topology coincide on $\bar{\mathcal{F}}$.

Le $G$ be a topological group.
Following \cite{article:DantzigMonotheticGroup}
(see also \cite[Definition 9.2 pag. 85]{book:HewittRossAHA_I})
we say that $G$ is a \emph{monothetic} group generated by $g$
if $g\in G$ and
the subgroup generated by $g$ is dense in $G$.
Of course if $G$ is monothetic generated by some element $g$ then $G$
is an abelian group.

We now recall a simple algebraic characterization of groups.
{
\def\sG{G}
\def\ea{g}
\def\eb{h}
\def\ec{k}
\def\en{e}
\def\em{f}
\def\ex{u}
\def\ey{v}
\begin{proposition}\label{thm::SemiToGroup}
Let $\sG$ be a not empty semigroup. Assume that for each $\ea,\eb\in\sG$
there exist $\ex,\ey\in\sG$ such that $\eb=\ex\ea=\ea\ey$.
Then $\sG$ is a group.
\end{proposition}

\proof
Let $\ea_0\in\sG$ be an arbitrarily chosen element of $\sG$.
Then we have $\ea_0=\en\ea_0$ and $\ea_0=\ea_0\em$ for some
$\en,\em\in\sG$.
We claim that for each $\ea\in\sG$ we have
$\ea=\en\ea$ and $\ea=\ea\em$.

Indeed, given  $\ea\in\sG$ there exists $\ex,\ey\in\sG$ such that
$\ea=\ea_0\ex=\ey\ea_0$, and hence
\begin{eqnarray*}
&&\en\ea=\en(\ea_0\ex)=(\en\ea_0)\ex=\ea_0\ex=\ea\\
&&\ea\em=(\ey\ea_0)\em=\ey(\ea_0\em)=\ey\ea_0=\ea.
\end{eqnarray*}
In particular we have $\en\em=\em$ and $\en\em=\en$, and hence $\en=\em$.

It follows that the element $\en$ is the unique element of $\sG$
which satisfies $\em\ea=\ea\em=\ea$ for each $\ea\in\sG$,
that is $\en$ is a neutral element for the semigroup $\sG$.

We end the proof showing that for each $\ea\in\sG$ there exists $\eb\in\sG$
such that $\ea\eb=\eb\ea=\en$.

Given $\ea\in\sG$ there exist $\eb,\ec\in\sG$ such that
$\eb\ea=\ea\ec=\en$. It suffices to prove that $\eb=\ec$.
Indeed we have
\begin{equation*}
\eb=\eb\en=\eb(\ea\ec)=(\eb\ea)\ec=\en\ec=\ec,
\end{equation*}
as desired.

\qed

}

Let $\pMap:\MSpace\to\MSpace$ be a not expansive map.

We denote by $\LimMapSet(\MSpace)$ (resp. $\LimMapSetEx(\MSpace)$)
the set of all continuous maps $\LMapA:\MSpace\to\MSpace$
(resp. $\LMapA:\MSpace\to\hat\MSpace$)
which are limit of a sequence of the iterates of the map $\pMap$,
that is there exist an increasing sequence of positive numbers
$\itx1<\itx2<\ldots$ such that the sequence
$$\pMap^{\itx1},\pMap^{\itx2},\ldots$$
converges uniformely on the compact subsets of $\MSpace$
to the map $\LMapA$.

We begin with the following easy lemma.

\begin{lemma}\label{stm:RecursetIsClosed}
Let $\MSpace$ be a metric space % with distance function $\distfunc_\MSpace$
and let $\pMap:\MSpace\to\MSpace$ be a not expansive map.
Then the set of $\pMap$-recurrent points of $\MSpace$ is closed in $\MSpace$.
\end{lemma}

\proof
Let $\distfunc_\MSpace$ be the distance function of $\MSpace$.

Let $\ppa_n\in \RecurrSet$ be a sequence of points converging to a point $\ppa\in\MSpace$.
We need to prove that then also $\ppa\in\RecurrSet$.
 
For each pair of positive integers $\nni,\nnj$ we have
\begin{eqnarray*}
\distfunc_\MSpace\bigl(\ppa,\pMap^\nnj(\ppa)\bigr)&\leq&
\distfunc_\MSpace(\ppa,\ppa_\nni)+
\distfunc_\MSpace\bigl(\ppa_\nni,\pMap^\nnj(\ppa_\nni)\bigr)+
\distfunc_\MSpace\bigl(\pMap^\nnj(\ppa_\nni),\pMap^\nnj(\ppa)\bigr)
\\
&\leq&
2\,\distfunc_\MSpace(\ppa,\ppa_\nni)+
\distfunc_\MSpace\bigl(\ppa_\nni,\pMap^\nnj(\ppa_\nni)\bigr).
\end{eqnarray*}

The quantities $\distfunc_\MSpace(\ppa,\ppa_\nni)$
and $\distfunc_\MSpace\bigl(\ppa_\nni,\pMap^\nnj(\ppa_\nni)\bigr)$
can be made arbitrarily small by suitable values of
$\nni$ and $\nnj$ with $\nnj$ arbitrarily large and hence
the quantity $\distfunc_\MSpace\bigl(\ppa,\pMap^\nnj(\ppa)\bigr)$
can be made arbitrarily small with
a suitable value of $\nnj$ arbitrarily large,
that is, by definition, $\ppa\in \RecurrSet$, as required.

\qed

The following proposition gives the basic properties
of the elements of $\LimMapSetEx(\MSpace)$.

\begin{proposition}\label{thm::LimStructure}
Let $\MSpace$ be a locally compact metric space
with countable basis and
let $\pMap:\MSpace\to\MSpace$ be a not expansive map.

Let $\CompDivSet$ be the set of the points $\ppa\in\MSpace$ having the
$\pMap$-orbit compactly divergent and
let $\RecurrSet$ be the set of $\pMap$-recurrent points of $\MSpace$.

Let $\pMapCA\in\LimMapSetEx(\MSpace)$ be given.

Then the following assertions hold:
\begin{enumerate}
\item $\pMapCA^{-1}(\infty)=\CompDivSet$;
\item $\pMapCA(\MSpace\setminus\CompDivSet)=\pMapCA(\RecurrSet)=\RecurrSet$.
\end{enumerate}
\end{proposition}

\proof
Let $\itx\itn$ be an increasing sequence of positive numbers such that
the sequece $\pMap^{\itx\itn}$ converges (pointwise) to $\pMapCA$.

Let $\ppa\in\CompDivSet$. Then the sequence
$\pMap^\kk(\ppa)$ is compactly divergent and hence
$\pMapCA(\ppa)=\limtoinfty\nu\pMap^{\itx\itn}(\ppa)=\infty$,
that is $\CompDivSet\sset\pMapCA^{-1}(\infty)$.
Conversely let $\ppa\in\pMapCA^{-1}(\infty)$.
Then the sequence $\pMap^\kk(x)$ is not relatively compact in
$\MSpace$ because $\limtoinfty\itn\pMap^{\itx\itn}(\ppa)=\infty$.
Theorem \ref{thm::OutOut} implies that the sequence $\pMap^\kk(\ppa)$
is compactly divergent, that is $\ppa\in\CompDivSet$.

This proves the first assertion of the Proposition.

Let us prove the second one.

From $\RecurrSet\sset\MSpace\setminus\CompDivSet$ it follows that
$\pMapCA(\RecurrSet)\sset\pMapCA(\MSpace\setminus\CompDivSet)$.

We end the proof showing that
$\pMapCA(\MSpace\setminus\CompDivSet)\subset \RecurrSet$
and
$\RecurrSet\subset\pMapCA(\RecurrSet)$.

Let $\ppa\in\MSpace\setminus\CompDivSet$.
It is not restrictive to assume that the sequence $\itx{\itn+1}-\itx\itn$
is increasing.
We have
\begin{equation*}
\limtoinfty\nu\pMap^{\itx\itn}(\ppa)=\pMapCA(\ppa)
\end{equation*}
and
\begin{equation*}
\limtoinfty\nu\pMap^{\itx{\itn+1}-\itx\itn}\bigl(\pMap^{\itx\itn}(\ppa)\bigr)=
\limtoinfty\nu\pMap^{\itx{\itn+1}}(\ppa)=\pMapCA(\ppa).
\end{equation*}
Proposition \ref*{thm::ConvLemma} then implies
\begin{equation*}
\limtoinfty\nu\pMap^{\itx{\itn+1}-\itx\itn}\bigl(\pMapCA(\ppa)\bigr)=\pMapCA(\ppa),
\end{equation*}
that is $\pMapCA(\ppa)\in \RecurrSet$.

Since $\ppa\in\MSpace\setminus\CompDivSet$ is arbitrary it follows that
$\pMapCA(\MSpace\setminus\CompDivSet)\subset \RecurrSet$.

Let now $\ppb\in\RecurrSet$.
We need to prove that there exists $\ppa\in \RecurrSet$ such that
$\pMapCA(\ppa)=\ppb$.

First observe that there exists an increasing sequence of positive integers
$\ity\itn$ such that
\begin{equation*}
\limtoinfty\nu\pMap^{\ity\itn}(\ppb)=\ppb.
\end{equation*}
We may assume that the sequence $\ity\itn-\itx\itn$ is increasing too.

We observe that the sequence $\pMap^\kk(\ppb)$ is not 
compactly divergent and hence, by Theorem \ref{thm::OutOut}, is
relatively compact in $\MSpace$.
Then we may assume that there exist $\ppa\in\MSpace$ such that
\begin{equation*}
\limtoinfty\nu\pMap^{\ity\itn-\itx\itn}(\ppb)=\ppa.
\end{equation*}
Lemma \ref{thm::CompToRecurrent} implies that $\ppa\in\RecurrSet$
and we also have
\begin{equation*}
\limtoinfty\nu\pMap^{\itx\itn}\bigl(\pMap^{\ity\itn-\itx\itn}(\ppb)\bigr)=
\limtoinfty\nu\pMap^{\ity\itn}(\ppb)=\ppb
\end{equation*}
and hence, by Proposition \ref{thm::ConvLemma},
\begin{equation*}
\pMapCA(\ppa)=\limtoinfty\nu\pMap^{\itx\itn}(\ppa)=\ppb.
\end{equation*}
as required.

\qed

As immediate consequence of the Proposition above is that
\begin{eqnarray*}
&&\LimMapSetEx(\MSpace\setminus\CompDivSet)=\LimMapSet(\MSpace\setminus\CompDivSet),\\
&&\LimMapSetEx(\RecurrSet)=\LimMapSet(\RecurrSet).
\end{eqnarray*}

It is straightforward to prove that $\LimMapSet(\MSpace\setminus\CompDivSet)$ and
$\LimMapSet(\RecurrSet)$ are commutative semigroups under the composition of maps.

The main result of this section is

\begin{theorem}\label{thm::GroupStructBase}
Let $\MSpace$ be a locally compact metric space with countable basis
and let $\pMap:\MSpace\to\MSpace$ be a not expansive map.
	
Let denote by $\RecurrSet$ and $\CompDivSet$ the set of points $\ppa\in\MSpace$
such that the $\pMap$-orbit of $\ppa$ is respectively relatively compact and
compactly divergent.

Then the following assertions hold:
{
\begin{enumerate}
\item\label{thm::GroupStructBase::itemBase}
	$\RecurrSet$ and $\CompDivSet$ are closed disjoint $\pMap$-invariant subset of $\MSpace$
	and the open set $\MSpace\setminus\CompDivSet$ also is $\pMap$-invariant;
	moreover we have $\RecurrSet\neq\void$ if, and only if, $\CompDivSet\neq\MSpace$;
\item\label{thm::GroupStructBase::itemNotEmpty}
	$\LimMapSetEx(\MSpace)$ is a not empty
	compact subset of $C(\MSpace,\hat{\MSpace})$;
\item\label{thm::GroupStructBase::itemDivergent}
	we have $\CompDivSet=\MSpace$ if, and only if,
	$\LimMapSetEx(\MSpace)$ contain the single map sending all $\MSpace$ to $\infty$.
\end{enumerate}
}
If $\CompDivSet\neq\MSpace$ then:
{
\renewcommand{\theenumi}{\roman{enumi}}
\begin{enumerate}
\item\label{thm::GroupStructBase::itemEx}
	the restriction map
	$$C(\MSpace,\hat{\MSpace})\ni\LMapA\mapsto\LMapA_{|\MSpace\setminus\CompDivSet}
	\in C(\MSpace\setminus\CompDivSet,\hat\MSpace)$$
	induces a homeomorphism between $\LimMapSetEx(\MSpace)$
	onto $\LimMapSet(\MSpace\setminus\CompDivSet)$;
\item\label{thm::GroupStructBase::itemGroupA}
	$\LimMapSet(\MSpace\setminus\CompDivSet)$ is a
	compact monothetic (abelian) topological group with respect to
	the composition of maps	generated by $\pMap\circ\rho$,
	where $\rho$ denotes the identity element of $\LimMapSet(\MSpace\setminus\CompDivSet)$;
\item\label{thm::GroupStructBase::itemSubgroup}
	the identity element $\rho\in\LimMapSet(\MSpace\setminus\CompDivSet)$ is a
	retraction of $\MSpace\setminus\CompDivSet$ onto the set of the recurrent points $\RecurrSet$;
\item\label{thm::GroupStructBase::itemIso}	
	the restriction map
	$$\LimMapSet(\MSpace\setminus\CompDivSet)\ni\LMapA\mapsto\LMapA_{|\RecurrSet}\in\LimMapSet(\RecurrSet)$$
	induces an isomorphism of topological groups between the group
	$\LimMapSet(\MSpace\setminus\CompDivSet)$ onto the group $\LimMapSet(\RecurrSet)$;
\item\label{thm::GroupStructBase::itemGroupB}
	$\LimMapSet(\RecurrSet)$ is a subgroup of the group of the surjective isometries
	of $\RecurrSet$ and is a compact monothetic topological group generated by
	the restriction of $\pMap$ to $\RecurrSet$. In particular the restriction
	of $\pMap$ to $\RecurrSet$ is a surjective isometry of $\RecurrSet$;
\item\label{thm::GroupStructBase::itemAction}
	the composition map
	$$\LimMapSet(\RecurrSet)\times\LimMapSet(\MSpace\setminus\CompDivSet)\ni
		(\LMapA,\LMapB)\mapsto\LMapA\circ\LMapB\in\LimMapSet(\MSpace\setminus\CompDivSet)
	$$
	is well defined and
	induces a (left) group action of $\LimMapSet(\RecurrSet)$ on
	$\LimMapSet(\MSpace\setminus\CompDivSet)$ which is free and transitive, that is
	$\LimMapSet(\MSpace\setminus\CompDivSet)$ is a principal homogeneous space for
	$\LimMapSet(\RecurrSet)$;
\item\label{thm::GroupStructBase::itemLimOrbit}
	for each $\ppa\in\MSpace\setminus\CompDivSet$ the set of all the
	accumulation points of the $\pMap$-orbit of $\ppa$ coincides with the
	$\LimMapSet(\RecurrSet)$-orbit of the point $\rho(\ppa)$.
\item\label{thm::GroupStructBase::itemLimUniqueOrbit}
	for each $\ppa,\ppb\in\MSpace\setminus\CompDivSet$ we have
	\begin{equation*}
		\limtoinfty\kk\pMap^{\kk}(\ppa)=\ppb
	\end{equation*}
	if, and only if, $\ppb=\rho(\ppa)$ and $\pMap(\ppb)=\ppb$.
\end{enumerate}
}
\end{theorem}

\proof
Let denote by $\distfunc_\MSpace$ the distance function on $\MSpace$.

(\ref{thm::GroupStructBase::itemBase}):
The invariance of the subseta $\RecurrSet$, $\CompDivSet$ and $\MSpace\setminus\CompDivSet$
is straightforward.

We already observed in lemma \ref{stm:RecursetIsClosed}
that $\RecurrSet$ is closed in $\MSpace$.
We prove that $\CompDivSet$ is closed showing that $\MSpace\setminus\CompDivSet$
is open.

Let $\ppa_0\in\MSpace\setminus\CompDivSet$.
By Theorem \ref{thm::OutOut} the orbit $\orbit$ of $\ppa_0$ is
relatively compact in $\MSpace$.
Choose $\e>0$ in such a way that $\orbit_\e$
is relatively compact in $\MSpace$.
Since the map $\pMap$ is not expansive it follows that
the orbit of each $\ppa\in\ball{\ppa_0}\e$ is contained in $\orbit_\e$
and hence $\ball{\ppa_0}\e\subset\MSpace\setminus\CompDivSet$.
Since $\ppa_0\in\MSpace\setminus\CompDivSet$ is arbitrary
it follows that $\MSpace\setminus\CompDivSet$ is open.

The last assertion follows directly from Lemma \ref{thm::CompToRecurrent}.

(\ref{thm::GroupStructBase::itemNotEmpty}): for each $n\in\NN$ set
\begin{equation*}
\mathcal{F}_n=\bigl\{\pMap^m\mid m\geq n\bigr\}.
\end{equation*}
Proposition \ref{thm::AbateComp} implies that $\bar{\mathcal{F}_n}$ is
a sequence of not empty compact subset of $C(\MSpace,\hat{\MSpace})$
such that $\bar{\mathcal{F}_n}\supset\bar{\mathcal{F}_{n+1}}\supset\cdots$
and
\begin{equation*}
\LimMapSetEx(\MSpace)=\bigcap_{n=1}^{\infty}\bar{\mathcal{F}_n}.
\end{equation*}
It follows that $\LimMapSetEx(\MSpace)$
is a non empty compact subset of $C(\MSpace,\hat{\MSpace})$
being the intersection of a decreasing family of not empty compact subset.

(\ref{thm::GroupStructBase::itemDivergent}):
assume that $\CompDivSet=\MSpace$.
Then for each $\ppa\in\MSpace$ the sequence $\ppa,\pMap(\ppa),\pMap^2(\ppa),\ldots$
is compactly divergent and hence if $\LMapA\in\LimMapSetEx(\MSpace)$
necessarily $\LMapA(x)=\infty$ for each $\ppa\in\MSpace$.

Conversely assume that the single map sending all $\MSpace$ to $\infty$
belongs to $\LimMapSetEx(\MSpace)$.
Then there exists an increasing sequence $\itx\itn$ of positive integers such that
for each $\ppa\in\MSpace$
\begin{equation*}
\limtoinfty\itn\pMap^{\itx\itn}(\ppa)=\infty.
\end{equation*}
It follows that the orbit of each point $\ppa\in\MSpace$ is not relatively compact
and hence, by Theorem \ref{thm::OutOut}, is compactly divergent, that is
$\CompDivSet=\MSpace$.

(\ref{thm::GroupStructBase::itemEx}):
Let $\LMapA\in\LimMapSetEx(\MSpace)$. Then Proposition \ref{thm::LimStructure}
implies that
\begin{equation*}
\LMapA(\MSpace\setminus\CompDivSet)=\RecurrSet\subset\MSpace\setminus\CompDivSet
\end{equation*}
and hence
$\LMapA_{\mid\MSpace\setminus\CompDivSet}\in\LimMapSet(\MSpace\setminus\CompDivSet)$.
Since $\LimMapSetEx(\MSpace)$ is compact and $\LimMapSet(\MSpace\setminus\CompDivSet)$
is Hausdorff it suffices to prove that the restriction map
$\LMapA\mapsto\LMapA_{|\MSpace\setminus\CompDivSet}$ (which is clearly continuous)
in injective and surjective.

Let $\LMapA,\LMapB\in\LimMapSetEx(\MSpace)$ and assume that
$\LMapA_{\mid\MSpace\setminus\CompDivSet}=\LMapB_{\mid\MSpace\setminus\CompDivSet}$.
Proposition \ref{thm::LimStructure} implies that for each $\ppa\in\CompDivSet$
\begin{equation*}
\LMapA(\ppa)=\LMapB(\ppa)=\infty
\end{equation*}
and hence $\LMapA=\LMapB$,
that is the restriction map 
$\LMapA\mapsto\LMapA_{|\MSpace\setminus\CompDivSet}$
is injective.

Let now $\LMapB\in\LimMapSet(\MSpace\setminus\CompDivSet)$
and let $\itx\itn$ be an increasing sequence of positive integers such that
for each $\ppa\in\MSpace\setminus\CompDivSet$
\begin{equation*}
\limtoinfty\itn\pMap^{\itx\itn}(\ppa)=\LMapB(\ppa).
\end{equation*}
Proposition \ref{thm::LimStructure} implies that for each $\ppa\in\CompDivSet$
\begin{equation*}
\limtoinfty\itn\pMap^{\itx\itn}(\ppa)=\infty
\end{equation*}
and hence Proposition \ref{thm::AbateComp} implies that the function
\begin{equation*}
\LMapA(x)=\left\{
\begin{array}{lll}
\LMapB(\ppa)&\ \ &\ppa\in\MSpace\setminus\CompDivSet\\
\infty&\ \ &\ppa\in\CompDivSet
\end{array}
\right.
\end{equation*}
belongs to $\LimMapSetEx(\MSpace)$ and clearly we have
$\LMapA_{|\MSpace\setminus\CompDivSet}=\LMapB$,
that is the restriction map 
$\LMapA\mapsto\LMapA_{|\MSpace\setminus\CompDivSet}$
is surjective.

(\ref{thm::GroupStructBase::itemGroupA}):
We already know that $\LimMapSet(\MSpace\setminus\CompDivSet)$
is a compact subset closed under the composition of maps.
We prove that $\LimMapSet(\MSpace\setminus\CompDivSet)$ is
a group using Proposition \ref{thm::SemiToGroup}.

Let $\LMapA,\LMapB\in\LimMapSet(\MSpace\setminus\CompDivSet)$.
It suffices to prove that there exists $\LMapC\in\LimMapSet(\MSpace\setminus\CompDivSet)$
such that $\LMapC\circ\LMapB=\LMapB\circ\LMapC=\LMapA$.

Let $\itx\itn$ and $\ity\itn$ be increasing sequences of positive integers such that
$\pMap^\itx\itn$ and $\pMap^\ity\itn$ converge respectively to
$\LMapA$ and $\LMapB$ on (the compact subset of) $\MSpace\setminus\CompDivSet$.

It is not restrictive to assume that the sequence $\itx\itn-\ity\itn$ is increasing
and the iterates $\pMap^{\itx\itn-\ity\itn}$
converge on $\MSpace\setminus\CompDivSet$
to a map $\LMapC\in\LimMapSet(\MSpace\setminus\CompDivSet)$.

For each $\ppa\in\MSpace\setminus\CompDivSet$ we have
\begin{eqnarray*}
&&\limtoinfty\itn\pMap^{\ity\itn}(\ppa)=\LMapB(\ppa),\\	
&&\limtoinfty\itn\pMap^{\itx\itn-\ity\itn}\bigl(\pMap^{\ity\itn}(\ppa)\bigr)
=\limtoinfty\itn\pMap^{\itx\itn}(\ppa)=\LMapA(\ppa)
\end{eqnarray*}
and hence Proposition \ref{thm::ConvLemma} implies that
\begin{eqnarray*}
&&\LMapC\bigl(\LMapB(\ppa)\bigr)=
\limtoinfty\itn\pMap^{\itx\itn-\ity\itn}\bigl(\LMapB(\ppa)\bigr)
=\LMapA(\ppa),
\end{eqnarray*}
that is $\LMapC\circ\LMapB=\LMapA$.
Since $\LimMapSet(\MSpace\setminus\CompDivSet)$ is commutative semigroup
we also have $\LMapB\circ\LMapC=\LMapA$, as required.

Let now $\itx\itn$ be an increasing sequence of positive integers such that the sequence
of the iterates $\pMap^{\itx\itn}$ converges to $\rho$, the unit element of
$\LimMapSet(\MSpace\setminus\CompDivSet)$.
Then the sequence $\pMap^{\itx\itn+1}$ converges to $\pMap\circ\rho$
and hence $\pMap\circ\rho\in\LimMapSet(\MSpace\setminus\CompDivSet)$.

Let $\LMapA\in\LimMapSet(\MSpace\setminus\CompDivSet)$ be arbitrary
and let $\ity\itn$ be an increasing sequence of positive integers such that
for each $\ppa\in\MSpace\setminus\CompDivSet$
\begin{eqnarray*}
&&\LMapA(\ppa)=\limtoinfty\itn\pMap^{\ity\itn}(\ppa).
\end{eqnarray*}
Then
\begin{eqnarray*}
&&\LMapA=\LMapA\circ\rho
=\limtoinfty\itn\pMap^{\ity\itn}\circ\rho
=\limtoinfty\itn(\pMap\circ\rho)^{\ity\itn}.
\end{eqnarray*}
Since $\LMapA\in\LimMapSet(\MSpace\setminus\CompDivSet)$ is arbitrary
it follows that $\LimMapSet(\MSpace\setminus\CompDivSet)$
is a monothetic group generated by $\pMap\circ\rho$.

(\ref{thm::GroupStructBase::itemSubgroup}):
let $\rho\in\LimMapSet(\MSpace\setminus\CompDivSet)$ be the identity element.
Then we have $\rho^2=\rho$ and hence $\rho$ is a retraction of
$\MSpace\setminus\CompDivSet$ onto its image, which by
Proposition \ref{thm::LimStructure} coincides with $\RecurrSet$.

(\ref{thm::GroupStructBase::itemIso}):
Let denote by $\varphi$ the restriction map
$$\LimMapSet(\MSpace\setminus\CompDivSet)\ni\LMapA
\mapsto\varphi(\LMapA)=\LMapA_{|\RecurrSet}\in\LimMapSet(\RecurrSet).$$
Of course $\varphi$ is a homomorphism of semigroups (with identity) between
the compact monothetic group $\LimMapSet(\MSpace\setminus\CompDivSet)$
and the Hausdorff topological semigroup $\LimMapSet(\RecurrSet)$.
It suffices then to prove that $\varphi$ is injective and onto.

Let us prove that $\varphi$ is injective. Since $\LimMapSet(\MSpace\setminus\CompDivSet)$
is a group it suffices to prove that the kernel of $\varphi$ is trivial.
Let $\LMapA\in\LimMapSet(\MSpace\setminus\CompDivSet)$
and assume that $\varphi(\LMapA)$ is the identity element of $\LimMapSet(\RecurrSet)$,
that is $\LMapA(\ppa)=\ppa$ for each $\ppa\in\RecurrSet$.

Let $\ppa\in\MSpace\setminus\CompDivSet$ and let
$\rho$ be the identity element of $\LimMapSet(\MSpace\setminus\CompDivSet)$.
Then $\rho(x)\in\RecurrSet$ and hence $\LMapA\bigl(\rho(\ppa)\bigr)=\rho(\ppa)$.
Since $\LMapA\circ\rho=\LMapA$ in $\LimMapSet(\MSpace\setminus\CompDivSet)$
we have
\begin{equation*}
\LMapA(\ppa)=\LMapA\bigl(\rho(\ppa)\bigr)=\rho(\ppa).
\end{equation*}
Since $\ppa\in\MSpace\setminus\CompDivSet$ is arbitrary it follows that $\LMapA=\rho$,
the identity element of $\LimMapSet(\MSpace\setminus\CompDivSet)$.

Let us prove that $\varphi$ is onto.
Let $\LMapB\in\LimMapSet(\RecurrSet)$ and let $\itx\itn$ be an increasing
sequence of positive integers such that $\pMap^{\itx\itn}$ converges
to $\LMapB$ on $\RecurrSet$.
Taking a subsequence if necessary we may assume that $\pMap^{\itx\itn}$ converges
to a map $\LMapA\in\LimMapSet(\MSpace\setminus\CompDivSet)$
which clearly satisfied $\varphi(\LMapA)=\LMapB$.

(\ref{thm::GroupStructBase::itemGroupB}):
Let $\rho$ be the identity element of the group $\LimMapSet(\MSpace\setminus\CompDivSet)$.
Since $\rho^2=\rho$ and by Proposition \ref{thm::LimStructure} also
$\rho(\MSpace\setminus\CompDivSet)=\RecurrSet$ it follows that
the identity element of the group $\LimMapSet(\RecurrSet)$, being the
restriction of $\rho$ to $\RecurrSet$, is the identity map of $\RecurrSet$.

Moreover $\pMap\circ\rho$ is a generator of $\LimMapSet(\MSpace\setminus\CompDivSet)$
and since the restriction of $\pMap\circ\rho$ to $\RecurrSet$ coincide with
the restriction of $\pMap$ to $\RecurrSet$ it follows that $\LimMapSet(\RecurrSet)$
is a compact monothetic group generated by $\pMap$.

Let now $\LMapA\in\LimMapSet(\RecurrSet)$ be arbitrary.
Since $\LimMapSet(\RecurrSet)$ is a group with unit element the
identity map of $\RecurrSet$ it follows
that $\LMapA^{-1}\in\LimMapSet(\RecurrSet)$ and hence
the image of $\RecurrSet$ under $\LMapA$ is all $\RecurrSet$.

Since $\LMapA$ and $\LMapA^{-1}$ are both not increasing for
each $\ppa,\ppb\in\RecurrSet$ we have
\begin{equation*}
\distfunc_\MSpace(\ppa,\ppb)\geq
\distfunc_\MSpace\bigl(\LMapA(\ppa),\LMapA(\ppb)\bigr)\geq
\distfunc_\MSpace\bigl(\LMapA^{-1}(\LMapA(\ppa)),\LMapA^{-1}(\LMapA(\ppb))\bigr)\geq
\distfunc_\MSpace(\ppa,\ppb),
\end{equation*}
and hence $\LMapA$ is an isometry of $\RecurrSet$ onto $\RecurrSet$.

(\ref{thm::GroupStructBase::itemAction}):
let $\LMapA\in\LimMapSet(\RecurrSet)$ and let
$\LMapB\in\LimMapSet(\MSpace\setminus\CompDivSet)$.
By Proposition \ref{thm::LimStructure} we have
$\LMapB(\MSpace\setminus\CompDivSet)=\RecurrSet$
and hence the composition $\LMapA\circ\LMapB$ is
well defined.

We need only to prove that given
$\LMapB,\LMapC\in\LimMapSet(\MSpace\setminus\CompDivSet)$
there exists a unique
$\LMapA\in\LimMapSet(\RecurrSet)$
such that $\LMapA\circ\LMapB=\LMapC$.

Let $\LMapB,\LMapC\in\LimMapSet(\MSpace\setminus\CompDivSet)$ be given.
If we choose $\LMapA$ as the restriction to $\LMapC\circ\LMapB^{-1}$ to $\RecurrSet$
(here $\LMapA^{-1}$ stands for the inverse of $\LMapA$ in the group
$\LimMapSet(\MSpace\setminus\CompDivSet)$) we clearly obtain that
$\LMapA\circ\LMapB=\LMapC$.

Let now $\LMapA_1,\LMapA_2\in\LimMapSet(\RecurrSet)$ and suppose
$\LMapA_1\circ\LMapB=\LMapA_2\circ\LMapB$.
Let $\ppb\in\RecurrSet$. By Proposition \ref{thm::LimStructure}
there exists $\ppa\in\MSpace\setminus\CompDivSet$
such that $\LMapB(\ppa)=\ppb$ and hence
\begin{equation*}
\LMapA_1(\ppb)
=\LMapA_1\bigl(\LMapB(\ppa)\bigr)
=\LMapA_2\bigl(\LMapB(\ppa)\bigr)
=\LMapA_2(\ppb).
\end{equation*}
Since $\ppb\in\RecurrSet$ is arbitrary then $\LMapA_1=\LMapA_2$.

(\ref{thm::GroupStructBase::itemLimOrbit}):
let $\ppa\in\MSpace\setminus\CompDivSet$ and let
$\ppb$ be an accumulation point of the $\pMap$-orbit of $\MSpace$.
Let $\itx\itn$ be an increasing sequence of positive integers
such that $\pMap^{\itx\itn}(\ppa)$ converges to $\ppb$.
Taking a subsequence if necessary we may suppose that the
sequence of functions $\pMap^{\itx\itn}$ converges to a
map $\LMapA\in\LimMapSet(\MSpace\setminus\CompDivSet)$.
Since $\LMapA\circ\rho=\LMapA$ in $\LimMapSet(\MSpace\setminus\CompDivSet)$
we have
\begin{equation*}
\ppb=\limtoinfty\itn\pMap^{\itx\itn}(\ppa)
    =\LMapA(\ppa)
    =\LMapA\bigl(\rho(\ppa)\bigr).
\end{equation*}
Since $\LMapA_{|\RecurrSet}\in\LimMapSet(\RecurrSet)$
and $\rho(\ppa)\in\RecurrSet$ it follows that
$\ppb$ is contained in the $\LimMapSet(\RecurrSet)$-orbit
of $\rho(\ppa)$.

Conversely assume that $\ppb=\LMapB\bigl(\rho(\ppa)\bigr)$
for some $\LMapB\in\LimMapSet(\RecurrSet)$.
Let $\LMapA\in\LimMapSet(\MSpace\setminus\CompDivSet)$ satisfying
$\LMapA_{|\RecurrSet}=\LMapB$ and let
$\itx\itn$ and $\itx\itm$ be
two increasing sequences of positive integers such that
$\pMap^{\itx\itn}$ and $\pMap^{\ity\itn}$ converges respectively
to $\LMapA$ and $\rho$.
Then
\begin{eqnarray*}
&&\limtoinfty\itn\pMap^{\ity\itn}(\ppa)=\rho(\ppa),\\
&&\limtoinfty\itn\pMap^{\itx\itn}(\rho\bigl(\ppa)\bigr)=
	\LMapA\bigl(\rho(\ppa)\bigr)=
	\LMapB\bigl(\rho(\ppa)\bigr)=\ppb,
\end{eqnarray*}
and hence, by Proposition \ref{thm::ConvLemma},
\begin{equation*}
\limtoinfty\itn\pMap^{\itx\itn+\ity\itn}(\ppa)
=\limtoinfty\itn\pMap^{\itx\itn}\bigl(\pMap^{\ity\itn}(\ppa)\bigr)
=\ppb
\end{equation*}
and this implies that $\ppb$ is an accumulation point of the
$\pMap$-orbit of $\ppa$.

(\ref{thm::GroupStructBase::itemLimUniqueOrbit}):
let $\ppa,\ppb\in\MSpace\setminus\CompDivSet$ and suppose
\begin{equation*}
	\limtoinfty\kk\pMap^{\kk}(\ppa)=\ppb.
\end{equation*}
Then $\ppb$ is the only accumulation point of the sequence
$\pMap^{\kk}(\ppa)$ and hence the $\LimMapSet(\RecurrSet)$-orbit
of $\rho(\ppa)$ contains the single element $\ppb$, that is
$\ppb=\rho(\ppa)$. Since the restriction of $\pMap$ to
$\RecurrSet$ belongs to $\LimMapSet(\RecurrSet)$ then
necessarily $\pMap(\ppb)=\ppb$.

Conversely assume $\ppb=\rho(\ppa)$ and $\pMap(\ppb)=\ppb$.
Since $\LimMapSet(\RecurrSet)$ is a monothetic group
generated by $\pMap$ it follows that the
$\LimMapSet(\RecurrSet)$-orbit of $\rho(\ppa)$
consists of the single element $\ppb=\rho(\ppa)$
and hence $\ppb$ is the only accumulation point of
the $\pMap$-orbit of $\ppa$ which by Theorem \ref{thm::OutOut}
is compact.
It is straightforward to prove that then the whole sequence
$\pMap^\kk(\ppa)$ converges to $\ppb$.

The proof of Theorem \ref{thm::GroupStructBase} is so completed.

\qed

\begin{remark} {\rm
We point out that the group $\LimMapSet(\MSpace\setminus\CompDivSet)$
is not a subgroup of the group of transformations of $\MSpace\setminus\CompDivSet$.
}
\end{remark}

\begin{remark} {\rm
$\LimMapSet(\RecurrSet)$ is a compact group of
surjective isometries of $\RecurrSet$.
We point out that in general the full group of the 
surjective isometries of a locally compact space in general
is not a locally compact topological group, unless the underlying space
is connected, as asserted by van Dantzig and an der Waerden Theorem \cite{article:DantzigWaerden};
see also \cite[Theorem 4.7 pag. 46]{book:KobayashiNomizuA}.
}
\end{remark}

\begin{remark} {\rm
$\LimMapSet(\MSpace\setminus\CompDivSet)$ and $\LimMapSet(\RecurrSet)$
are (compact) monothetic topological group.
It is a standard result of abstract harmonic analysis that
a locally compact monothetic topological group either is 
isomorphic to $\ZZ$ or is compact; see, e. g.,
\cite[Theorem 9.1, pag 84]{book:HewittRossAHA_I}.
}
\end{remark}

\begin{remark} {\rm
The proof that $\LimMapSet(\MSpace\setminus\CompDivSet)$
is a group given above % in \ref{thm::GroupStructBase::itemGroupA}
is a simplified adaptation of the proof of a more general result 
on the existence of groups in compact topological semigroups
given in \cite{article:NumakuraBicompactGroups};
see also \cite[Lemma 9.17, pag 100]{book:HewittRossAHA_I}.
}
\end{remark}

In analogy with \cite[pag. 145]{book:AbateBook} we say that the
unit element $\rho$ in the group $\LimMapSet(\MSpace\setminus\CompDivSet)$
is the \emph{limit retraction} of $\pMap$ and
we define the \emph{extended limit retraction} of $\pMap$ as the unique map
$\hat{\rho}\in\LimMapSetEx(\MSpace)$ which coincides with $\rho$ on
$\MSpace\setminus\CompDivSet$.

We then have the following characterization of $\rho$. % and $\hat{\rho}$.

\begin{proposition}
Let $\MSpace$ be a locally compact metric space with countable basis
and let $\pMap:\MSpace\to\MSpace$ be a not expansive map.

Assume that $\MSpace$ contains at least a $\pMap$-recurrent point.
Then the retraction limit
of $\pMap$ is the unique element of 
$\LimMapSet(\MSpace\setminus\CompDivSet)$
which leaves invariant
each $\pMap$-recurrent point of $\MSpace$.
\end{proposition}

\proof
Let $\rho$ be the retraction limit of $\pMap$.
Then $\rho(\ppa)=\ppa$ if $\ppa$ is a $\pMap$-recurrent point of $\MSpace$.

Conversely, let $\LMapA\in\LimMapSet(\MSpace\setminus\CompDivSet)$
be a map which leaves invariant each $\pMap$-recurrent point of $\MSpace$.
Then $\rho(\ppa)$ is $\pMap$-recurrent for each 
$\ppa\in\MSpace\setminus\CompDivSet$ and hence
\begin{equation*}
\LMapA(\ppa)=
(\LMapA\circ\rho)(\ppa)=
\LMapA\bigl(\rho(\ppa)\bigr)=\rho(\ppa),
\end{equation*}
that is $\LMapA=\rho$.

\qed

\begin{theorem}\label{thm::Convergence}
Let $\MSpace$ be a locally compact metric space
with countable basis and
let $\pMap:\MSpace\to\MSpace$ be a not expansive map.

Then the sequence $\pMap^\kk$ converges on $\MSpace$
to the extended retraction map $\hat{\rho}$
if, and only if,
$\pMap(\ppa)=\ppa$ for each $\pMap$-recurrent point
of $\MSpace$.
\end{theorem}

\proof
Let $\CompDivSet=\hat{\rho}^{-1}(\infty)$ and let
$\RecurrSet$ be the set of all the $\pMap$-recurrent 
points of $\MSpace$.

Assume that $\pMap^\kk$ converges to $\hat{\rho}$
and let $\ppb\in\RecurrSet$ be arbitrary.
Then $\hat{\rho}(\ppb)=\ppb$ and
by hypothesis $\pMap^{\kk}(\ppb)$ converges to $\hat{\rho}(\ppb)$,
that is
\begin{equation*}
\limtoinfty\kk\pMap^\kk(\ppb)=\hat{\rho}(\ppb)=\ppb
\end{equation*}
and hence
\begin{equation*}
\pMap(\ppb)=\limtoinfty\kk\pMap^{\kk+1}(\ppb)=\ppb.
\end{equation*}

Conversely assume that $\pMap(\ppb)=\ppb$ for each $\ppb\in\RecurrSet$
and let $\ppa\in\MSpace$ be arbitrary.
If $\ppa\in\CompDivSet$ then 
\begin{equation*}
\limtoinfty\kk\pMap^\kk(\ppa)=\infty=\hat{\rho}(\ppa).
\end{equation*}
Suppose now that $\ppa\in\MSpace\setminus\CompDivSet$
and set $\ppb=\hat{\rho}(\ppa)=\rho(\ppa)$.
Then by Proposition \ref{thm::LimStructure} we have $\ppb\in\RecurrSet$
and hence $\pMap(\ppb)=\ppb$.

The statement (\ref{thm::GroupStructBase::itemLimUniqueOrbit})
of Theorem \ref{thm::GroupStructBase} implies that
\begin{equation*}
\limtoinfty\kk\pMap^\kk(\ppa)=\ppb=\hat{\rho}(\ppa)
\end{equation*}
and we are done.

\qed

We end this section with a simple consequence of the
(existence of the)
extended retraction.

\begin{proposition}
Let $\MSpace$ be a locally compact metric space
with countable basis and
let $\pMap:\MSpace\to\MSpace$ be a not expansive map.

Let $\CompDivSet$ be the set of the points $\ppa\in\MSpace$ having the
$\pMap$-orbit compactly divergent and
let $\RecurrSet$ be the set of $\pMap$-recurrent points of $\MSpace$.

If $\RecurrSet$ is compact then $\CompDivSet$ is open and closed in $\MSpace$.

In particular if $\MSpace$ is (not empty and) connected and
$\RecurrSet$ is compact then either $\CompDivSet=\MSpace$ and
$\RecurrSet=\void$ or $\CompDivSet=\void$ and $\RecurrSet$ is
not empty and connected.
\end{proposition}

\proof
We already know that $\CompDivSet=\hat{\rho}^{-1}(\infty)$ is closed.
If $\RecurrSet$ is compact then $\hat{\MSpace}\setminus\RecurrSet$
is open in $\hat{\MSpace}$ and hence
$\CompDivSet=\hat{\rho}^{-1}(\hat{\MSpace}\setminus\RecurrSet)$
is also open.

\qed

\section{\label{section::ComplexSpaces}Complex hyperbolic spaces}
\def\pcra#1#2{\left|\dfrac{\displaystyle{{#1}-{#2}}}{\displaystyle{1-\bar{#1}{#2}}}\right|}
We now recall some basic fact on the theory of
Kobayashi hyperbolic complex spaces.
For more details and further results see, e. g.
\cite{book:Kobayashi} or \cite{book:ComplexHyperbolicSpaces}.

Let $$\Delta=\bigl\{z\in\CC\mid\abs{z}<1\bigr\}.$$
The \emph{Poicar\'e metric} on $\Delta$ is the Riemannian metric given by
\begin{equation*}
\dfrac{dzd\bar z}{\bigl(1-\abs{z}^2\bigr)^2}
\end{equation*}
and the associated distance is given by
\begin{equation*}
\omega(z,w)=\frac{1}{2}
\log\frac{1+\pcra{z}{w}}{1-\pcra{z}{w}},\ z,w\in\Delta
\end{equation*}

It is well-known that each holomorhic map $\pMap:\Delta\to\Delta$ is
not increasing, that is for each $z,w\in\Delta$ we have
\begin{equation*}
\omega\bigl(\pMap(z),\pMap(w)\bigr)\leq\omega(z,w).
\end{equation*}

Let now $\CSpace$ be a connected complex space.
An \emph{analytic chain}
$$\alpha=\{z_0,\ldots,z_m;w_0,\ldots,w_m;\varphi_0,\ldots,\varphi_m\}$$
connecting two points $\ppa$ and $\ppb$ of $\CSpace$ is
a sequence of points $z_0,\ldots,z_m,w_0,\ldots,w_m\in\Delta$
and holomorphic maps $\varphi_0,\ldots,\varphi_m:\Delta\to\CSpace$
such that
$\varphi_0(z_0)=\ppa$,
$\varphi_j(w_j)=\varphi_{j+1}(z_{j+1})$ for $j=0,\ldots,m-1$ and
$\varphi_m(w_m)=\ppb$.
The \emph{length} of the chain $\alpha$ is
\begin{equation*}
\omega(\alpha)=\sum_{j=0}^m\omega(z_j,w_j).
\end{equation*}
The \emph{Kobayashi (pseudo)-distance} $k_\CSpace(\ppa,\ppb)$ between
the two points $\ppa$ and $\ppb$ is the infimum of the lengths of
the analytic chains connecting the points $\ppa$ and $\ppb$.

The complex space $\CSpace$ is \emph{hyperbolic} (in the sense of Kobayashi)
if $k_\CSpace(\ppa,\ppb)>0$ for each pair of distinct points $\ppa,\ppb\in\CSpace$.
In this case $k_\CSpace$ is a distance function on $\CSpace$ which induces on
$\CSpace$ its original topology.

Moreover any holomorphic map $\hMap:\CSpace\to\CSpace$ is not expansive
with respect to the Kobayashi distance, that is, if $\ppa$ and $\ppb$ are
points of $\CSpace$ then
\begin{equation*}
k_\CSpace\bigl(\pMap(\ppa),\pMap(\ppb)\bigr)\leq k_\CSpace(\ppa,\ppb).
\end{equation*}

It is clear that Theorem \ref{thm::HyperbolicSpaces}
is an immediate consequence of
Theorem \ref{thm::OutOut}.

We end this paper with the following:

\begin{theorem}
Let $\CSpace$ be a connected hyperbolic complex space
and let $\pMap:\CSpace\to\CSpace$ be an holomorphic map.

Then the set $\RecurrSet$ of the $\pMap$-recurrent points of $\CSpace$
is a closed complex subspace of $\CSpace$ and
each singular point of $\RecurrSet$ also is a singular point of $\CSpace$.
\end{theorem}

\proof
If $\RecurrSet=\void$ there is nothing to prove.
So assume that $\RecurrSet\neq\void$.
Let $\CompDivSet$ be the set of points $\ppa\in\MSpace$
such that the $\pMap$-orbit of $\ppa$ compactly divergent.
We already know that $\RecurrSet$ and $\CompDivSet$ are
closed disjoint subset of $\CSpace$.

Let $\rho:\CSpace\setminus\CompDivSet\to\RecurrSet$ be
the limit retraction of $\pMap$.
The map $\rho$ is holomorphic being the limit of a sequence
of holomorphic functions.
The set $\CSpace\setminus\CompDivSet$ is open in $\CSpace$ and
\begin{equation*}
\RecurrSet=\bigl\{\ppa\in\CSpace\setminus\CompDivSet\mid
	\rho(\ppa)=\ppa\bigr\}
\end{equation*}
and hence $\RecurrSet$ is a complex subspace of 
the open set $\CSpace\setminus\CompDivSet$.

Since $\RecurrSet$ is closed in $\CSpace$ and contained
in the open set $\CSpace\setminus\CompDivSet$
it is then a complex subspace of $\CSpace$.

We end the proof showing that each point of $\RecurrSet$
which is regular point of $\CSpace$ also is a
regular point of $\RecurrSet$.

Let $\ppa\in\RecurrSet$ and assume that $\ppa$ is
a regular point of $\CSpace$.
Then there is a suitable connected neighbourhood $\ngh$ of $\ppa$
in $\CSpace$ such that each point of $\ngh$ is 
a regular point of $\CSpace$, that is $\ngh$ is a complex manifold,
and $\rho(\ngh)=\RecurrSet\cap\ngh$.

By a result of Rossi (\cite[Theorem 7.1]{article:HugoRossiSmoothRetract})
it follows that $\RecurrSet\cap\ngh$ is a smooth sub-manifold of $\ngh$,
that is $\ppa$ is a regular point of $\RecurrSet$.

\qed

\bibliographystyle{amsalpha}
\providecommand{\bysame}{\leavevmode\hbox to3em{\hrulefill}\thinspace}
\providecommand{\MR}{\relax\ifhmode\unskip\space\fi MR }
\providecommand{\MRhref}[2]{%
  \href{http://www.ams.org/mathscinet-getitem?mr=#1}{#2}
}
\providecommand{\href}[2]{#2}

\end{document}